\documentclass[conference]{IEEEtran}
\IEEEoverridecommandlockouts
\usepackage{cite}
\usepackage{amsmath,amsthm,amssymb,bbm,bm,mathtools}
\usepackage{subcaption}
\usepackage{graphicx}
\usepackage{textcomp}
\usepackage{xcolor}

\newtheorem{theorem}{Theorem}
\newtheorem{lemma}{Lemma}
\newtheorem{remark}{Remark}

\newcommand{\trace}{\mathtt{tr}}
\newcommand{\Expectation}{\mathbb{E}}
\renewcommand{\Re}{\mathbb{R}}

\def\BibTeX{{\rm B\kern-.05em{\sc i\kern-.025em b}\kern-.08em
    T\kern-.1667em\lower.7ex\hbox{E}\kern-.125emX}}

\title{Computationally Efficient Covariance Steering for Systems Subject to Parametric Disturbances and Chance Constraints}

\author{Jacob W. Knaup and Panagiotis Tsiotras, \IEEEmembership{Fellow, IEEE}
\thanks{J. Knaup is with the School of Interactive Computing, College of Computing and
the Institute for Robotics and Intelligent Machines,
Georgia Institute of Technology,  Atlanta, GA 30332--0250 USA (e-mail: jacobk@gatech.edu)}
\thanks{P. Tsiotras is with the School of Aerospace Engineering and the
 Institute for Robotics and Intelligent Machines,
Georgia Institute of Technology,  Atlanta, GA 30332--0150 USA (e-mail: tsiotras@gatech.edu)}
}

\begin{document}

\maketitle

\begin{abstract}
This work investigates the finite-horizon optimal covariance steering problem for discrete-time linear systems subject to both additive and multiplicative uncertainties as well as state and input chance constraints. In particular, a tractable convex approximation of the optimal covariance steering problem is developed by tightening the chance constraints and by introducing a suitable change of variables. 
The solution of the convex approximation is shown to be a valid (albeit potentially suboptimal) solution to the original chance-constrained covariance steering problem.
\end{abstract}


\section{Introduction}

Covariance control examines the problem of driving a stochastic system from a given initial distribution to a specified final distribution. The problem has been previously studied extensively for the infinite-horizon case \cite{hotz1987covariance, iwasaki1992quadratic, xu1992improved}, but has only recently been studied for the finite-horizon case (referred to as covariance steering). Specifically, in \cite{goldshtein2017finite, chen2015optimal1, chen2015optimal2, chen2018optimal}, the authors introduced the finite-horizon covariance steering problem. In \cite{bakolas2016optimal, bakolas2018finite}, the authors added expectation constraints, reference \cite{okamoto2018optimal} introduced state chance constraints, and in \cite{okamoto2019input} the authors considered hard input constraints. Additionally, covariance steering has been applied to robotics path planning in \cite{okamoto2019optimal, zheng2022belief}, differential games in \cite{chen2019covariance}, and re-entry, descent, and landing tasks for space operations in \cite{ridderhof2018uncertainty, goyal2021optimal, ridderhof2022chance}.
However, most recent work on covariance steering focuses on linear systems subject only to additive noise, where the initial and final distributions as well as the noise characteristics are given as Gaussian distributions. 
A few of notable exceptions are \cite{liu2022generic} which considered additive generic (non-Gaussian) noise and \cite{ridderhof2019nonlinear, chen2021covariance} which considered covariance control for nonlinear systems.

We examine the problem of steering a stochastic linear system from an initial distribution characterized by its first two moments to a given terminal mean and covariance in finite-time, when the system is subject to parametric uncertainties (i.e. the disturbances enter both multiplicatively with the state and control as well as additively). This work is most similar to \cite{liu2022optimal}, except that the authors of \cite{liu2022optimal} derive an algebraic Riccati equation-based solution for the continuous time case and do not consider chance constraints. In contrast, this work aims at computationally efficient solutions to the chance-constrained discrete-time case, which has only previously been addressed for systems with purely additive (not multiplicative) disturbances.

Although the problem of steering the first two moments of a linear system subject to state and control chance constraints and mixed additive and multiplicative noise from initial to given final conditions is, 
in general, nonlinear, this work shows that the problem may be represented by a tightened convex problem formulation, the solution of which ensures the satisfaction of the original problem's chance and terminal constraints. Therefore, the optimal solution of the convex reformulation is an admissible (albeit potentially sub-optimal) solution of the original nonlinear program. The convex reformulation is accomplished by employing Boole's and Cantelli's inequalities, and by relaxing the covariance propagation to a semi-definite constraint, similar to the approach used in \cite{farina2016model} for deriving an approximate convex problem for a model predictive control scheme. 

During the preparation of this manuscript, the authors became aware of a similar study \cite{balci2022covariance}. 
The authors of \cite{balci2022covariance} also studied the covariance steering problem for linear systems affected by multiplicative disturbances and solved the problem using semi-definite programming. 
However, whereas the work of \cite{balci2022covariance} assumes  multiplicative disturbances affecting the state and control, as well as the additive disturbances, they 
are all independent, which simplifies the covariance propagation.
Contrary, the current work does not assume that the disturbances affecting the system at a given time step are independent of each other. 
Additionally, the proposed approach in this paper includes chance constraints in the problem formulation and shows how they can be ensured using Boole's and Cantelli's inequalities and incorporated into the semi-definite program using a linear bounding technique, whereas the work of \cite{balci2022covariance} does not address chance constraints. 

The rest of this paper is organized as follows. 
In Section~\ref{sec:problem_formulation}, the optimal covariance steering problem is presented for linear systems subject to additive and multiplicative disturbances as well as state and control chance constraints. 
Section~\ref{sec:covariance_steering_controller_design} contains the main result of this work which is an admissible convex 
reformulation of the optimal covariance steering problem, accomplished by arriving at a deterministic problem in Section \ref{sec:deterministic_reformulation} and by convexifying the problem via a relaxation of 
the covariance propagation equations and deriving linear bounds on the probability of constraint violations in Section~\ref{sec:convex_reformulation}. 
Numerical simulations verifying the proposed approach are presented in Section~\ref{sec:numerical_results}, and Section~\ref{sec:conclusion} presents closing remarks. 

\subsubsection*{Notation}

This paper employs several standard notation practices. A random variable drawn from a normal distribution with mean $\mu$ and covariance matrix $\Sigma$ is denoted by $x \sim \mathcal{N}(\mu, \Sigma)$. $\Expectation[\cdot]$ denotes the expectation operator, and $\Pr(x)$ denotes the probability of event $x$. $I_n$ denotes the $n \times n$ identity matrix, and $\mathtt{tr}(\cdot)$ denotes the trace operation. A symmetric positive (semi)-definite matrix is denoted by $M \succ 0$ ($M \succeq 0$).

\section{Problem Formulation} \label{sec:problem_formulation}

Consider the system 
\begin{align} \label{eq:stochastic_sys}
    x_{k+1} &= A_k x_k + B_k u_k + d_k,
\end{align}
where $x_k \in \mathbb{R}^{n_x}$, $u_k \in \mathbb{R}^{n_u}$. 
Let the initial conditions be given as $\Expectation[x_0] = \mu_I$ and $\Expectation[(x_0 - \Expectation[x_0])(x_0 - \Expectation[x_0])^\top] = \Sigma_I$, where $\Sigma_I \succ 0$.
Additionally, the system matrices are comprised of a constant, known component, and a time-varying stochastic component given by
\begin{subequations}
\begin{align}
    A_k &= \bar{A} + \sum_{j=1}^{m} \tilde{A}_{j} q_{j, k}, \\
    B_k &= \bar{B} + \sum_{j=1}^{m} \tilde{B}_{j} q_{j, k}, \\
    d_k &= \bar{d} + \sum_{j=1}^{m} \tilde{d}_{j} q_{j, k},
\end{align}
\end{subequations}
where $q_{j, k} \sim \mathcal{N}(0, 1)$, for all $k = 0, 1, \dots$, is a Gaussian white noise. 
Therefore,
\begin{subequations}\label{eq:noise_independence}
\begin{align}
    \mathbb{E}[q_{j_1, k_1} q_{j_2, k_2}] = 
    \begin{cases}1  & \text{where } j_1 = j_2 \text{ and } k_1 = k_2, \\
    0 & \text{otherwise}.
    \end{cases}
\end{align}
Furthermore, we assume
\begin{align}
    \mathbb{E}[x_k q_{j, k}] = \mathbb{E}[x_k] \Expectation[q_{j, k}] = 0, \\
    \mathbb{E}[u_k q_{j, k}] = \Expectation[u_k] \Expectation[q_{j, k}] = 0,
\end{align}
\end{subequations}
for $k = 0, 1, \dots, j = 1, \dots, m$ which stem from causality considerations. 

The state and control inputs in (\ref{eq:stochastic_sys}) are subject to the chance constraints
\begin{align}\label{eq:state_and_input_const}
	\Pr(x_k \in \mathcal{X}) \geq 1 - p_x,\quad \Pr(u_k \in \mathcal{U}) \geq 1 - p_u,
\end{align}
for all $k = 0, 1, \dots, N-1$, where $\mathcal{X} \subseteq \Re^{n_x}$ and $\mathcal{U}\subseteq \Re^{n_u}$ are convex sets and $p_x, p_u \in (0, 0.5)$ are given maximal probabilities of constraint violation.
We further assume that the sets $\mathcal{X}$ and $\mathcal{U}$ can be written as an intersection of a finite number of linear inequality constraints as follows
\begin{align}
	\mathcal{X} &\triangleq 
	\bigcap_{i_x = 1}^{N_s} \left\{x: \alpha_{x,i_x}^\top x \leq \beta_{x,i_x}\right\},\label{eq:Xdefinition}\\
	\mathcal{U} &\triangleq 
	\bigcap_{i_u = 1}^{N_c} \left\{u: \alpha_{u,i_u}^\top u \leq \beta_{u,i_u}\right\},\label{eq:Udefinition}
\end{align}
where $\alpha_{x,i_x}\in\Re^{n_x}$ and $\alpha_{u,i_u}\in\Re^{n_u}$ are constant vectors, and $\beta_{x,i_x} \geq 0$ and $\beta_{u,i_u} \geq 0$ are constant scalars.

We wish to steer \eqref{eq:stochastic_sys} to a given final mean $\mu_F \in \mathcal{X}$ and covariance $\Sigma_F \succ 0$ at time $N$, such that 
\begin{align*}
    \Expectation[x_N] = \mu_F, \quad \Expectation[(x_N - \Expectation[x_N])(x_N - \Expectation[x_N])^\top] = \Sigma_F,
\end{align*}
while minimizing the cost function
\begin{align}\label{eq:cost_function}
    J(\mu_I, \Sigma_I; u_0, \dots, u_{N-1}) &= \Expectation\left[\sum_{k=0}^{N-1} \ell(x_k, u_k) \right].
\end{align}
In particular, we will examine the case where $\ell(\cdot, \cdot)$ has the special form
\begin{align}\label{eq:stage_cost}
    \ell(x, u) &= x^\top Q x + u^\top R u.
\end{align}
The problem may thus be summarized as follows.
Given $\mu_I, \Sigma_I, \mu_F, \Sigma_F$, determine the control sequence $\mathbf{u} = \{u_0, \dots, u_{N-1}\}$ that solves the following finite-time, optimal covariance steering problem
\begin{subequations}\label{prob:cs_prob}
\begin{align}
    &\min_\mathbf{u} ~ J(\mu_I, \Sigma_I; \mathbf{u}) = \Expectation\left[\sum_{k=0}^{N-1} x_k^\top Q x_k + u_k^\top R u_k \right], \label{eq:cs_cost}\\
    \text{s}&\text{ubject to } \nonumber\\
    &\Expectation[x_0] = \mu_I, \label{eq:cs_init_mean}\\
    &\Expectation[(x_0 - \Expectation[x_0])(x_0 - \Expectation[x_0])^\top] = \Sigma_I, \label{eq:cs_init_cov}\\
    &x_{k+1} = (\bar{A} + \sum_{j=1}^{m} \tilde{A}_{j} q_{j, k}) x_k
    + (\bar{B} + \sum_{j=1}^{m} \tilde{B}_{j} q_{j, k}) u_k \nonumber\\
    &+ \bar{d} + \sum_{j=1}^{m} \tilde{d}_{j} q_{j, k}, \quad k = 0, \dots, N-1 \label{eq:cs_dynamics}\\
    &\Pr(x_k \in \mathcal{X}) \geq 1 - p_x, \quad k = 0, \dots, N-1 \label{eq:cs_state_chance_const}\\
    &\Pr(u_k \in \mathcal{U}) \geq 1 - p_u, \quad k = 0, \dots, N-1 \label{eq:cs_control_chance_const}\\
    &\Expectation[x_N] = \mu_F, \label{eq:cs_terminal_mean_const}\\
    &\Expectation[(x_N - \Expectation[x_N])(x_N - \Expectation[x_N])^\top] = \Sigma_F. \label{eq:cs_terminal_cov_const}
\end{align}
\end{subequations}

\section{Covariance Steering Controller Design} \label{sec:covariance_steering_controller_design}

We now turn out attention to formulating a computationally tractable approximation of Problem \eqref{prob:cs_prob} using standard relaxations, the solution of which will provide a (suboptimal) feasible solution to the original problem. 
Specifically, we first formulate Problem \eqref{prob:cs_prob} as a deterministic optimal control problem by evaluating the expectations in \eqref{eq:cs_cost}-\eqref{eq:cs_terminal_cov_const} and derive explicit expressions for the propagation of the system mean and covariance. 
Additionally, we use Boole's inequality and Cantelli's inequality to approximate the chance constraints \eqref{eq:cs_state_chance_const}-\eqref{eq:cs_control_chance_const} as deterministic inequality constraints. We then show that the deterministic problem can be cast as a convex programming problem by performing a change of variables, tightening the chance constraints, and relaxing the terminal covariance constraint \eqref{eq:cs_terminal_cov_const} to a linear matrix inequality (LMI) constraint.

\subsection{Expectation and Uncertainty Propagation} \label{sec:uncertainty_propagation}

We may write the nominal system as 
\begin{align}
    \mathbb{E}\left[x_{k+1}\right] &= \mathbb{E}\left[A_k x_k + B_k u_k + d_k \right] \nonumber\\
    &= \mathbb{E}[(\bar{A} + \sum_{j=1}^{m} \tilde{A}_{j} q_{j, k} ) x_k 
    + (\bar{B} + \sum_{j=1}^{m} \tilde{B}_{j} q_{j, k} ) u_k \nonumber\\
    &+ \bar{d} + \sum_{j=1}^{m} \tilde{d}_{j} q_{j, k} ].
\end{align}
Using the independence of the disturbances \eqref{eq:noise_independence} and the fact that $\mathbb{E}[q_{j, k}] = 0$ for all $k = 0, 1, \dots$ and $j = 1, \dots, m$, the nominal system is given by
\begin{align}\label{eq:mean_dynamics}
    \bar{x}_{k+1} &= \bar{A} \bar{x}_k + \bar{B} \bar{u}_k + \bar{d},
\end{align}
where $\bar{x}_k = \mathbb{E}[x_k]$ and $\bar{u}_k = \mathbb{E}[u_k]$.

Next, let us define the deviation of the stochastic system from the nominal dynamics as $\tilde{x}_k = x_k - \bar{x}_k$. The error dynamics are given as 
\begin{align}\label{eq:x_tilde_dynamics}
    \tilde{x}_{k+1} &= 
    \bar{A} \tilde{x}_k + \bar{B} \tilde{u}_k \nonumber\\
    &
    + \sum_{j=1}^{m} \tilde{A}_{j} q_{j, k} x_k 
    + \sum_{j=1}^{m} \tilde{B}_{j} q_{j, k} u_k
    + \sum_{j=1}^{m} \tilde{d}_{j} q_{j, k},
\end{align}
where $\tilde{u}_k = u_k - \bar{u}_k$.
Considering \eqref{eq:noise_independence}, the covariance dynamics are given by
\begin{align} \label{eq:covariance_dynamics}
    &\Sigma_{x_{k+1}} = \bar{A} \Sigma_{x_k} \bar{A}^\top + \bar{A} \Sigma_{x_k u_k} \bar{B}^\top + \bar{B} \Sigma_{x_k u_k}^\top \bar{A}^\top + \bar{B} \Sigma_{u_k} \bar{B}^\top \nonumber\\
    &\hspace{-2mm}+ \sum_{j=1}^{m} (\tilde{A}_j \Sigma_{x_k} \tilde{A}_j^\top
    + \tilde{A}_j \Sigma_{x_k u_k} \tilde{B}_j^\top 
    + \tilde{B}_{j} \Sigma_{x_k u_k}^\top \tilde{A}_j^\top
    + \tilde{B}_{j} \Sigma_{u_k} \tilde{B}_j^\top) \nonumber\\
    &+ \sum_{j=1}^{m} (\tilde{A}_j \bar{x}_k + \tilde{B}_j \bar{u}_k + \tilde{d}_j)(\tilde{A}_j \bar{x}_k + \tilde{B}_j \bar{u}_k + \tilde{d}_j)^\top
\end{align}
where $\Sigma_{u_k} = \Expectation[(u_k - \Expectation[u_k])(u_k - \Expectation[u_k])^\top]$, and $\Sigma_{x_k u_k} = \Expectation[(x_k - \Expectation[x_k])(u_k - \Expectation[u_k])^\top]$ are properties of the particular control policy under consideration. 

\subsection{Control Policy} \label{sec:control_policy}

Rather than optimizing over control actions, we optimize over control policies.
In order to steer the mean and covariance of the system dynamics, we propose the following affine state-feedback control policy 
\begin{align} \label{eq:state_feedback}
u_k = L_k \tilde{x}_k + c_k,
\end{align}
where $L_k \in \Re^{n_u \times n_x}$ and $c_k \in \Re^{n_u}$ are new optimization variables, for $k = 0, \dots, N-1$.
Thus, 
\begin{subequations}
\begin{align}
    \bar{u}_k &= c_k ,\\
    \Sigma_{u_k} &= L_k \Sigma_{x_k} L_k^\top ,\\
    \Sigma_{x_k u_k} &= \Sigma_{x_k} L_k^\top,
\end{align}
\end{subequations}
for $k = 0, \dots, N-1$.

\subsection{Deterministic Reformulation} \label{sec:deterministic_reformulation}

We are now ready to formulate a deterministic approximation of Problem \eqref{prob:cs_prob}.
First, the expectation of the stage cost function \eqref{eq:stage_cost} can be written in terms of the mean and covariance as 
    $\Expectation[\ell(x_k, u_k)] 
    = \bar{x}_k^\top Q \bar{x}_k + \trace(Q \Sigma_{x_k}) + \bar{u}_k^\top R \bar{u}_k + \trace(R \Sigma_{u_k})$.
Next, using Boole's inequality \cite{prekopa1988boole} and \eqref{eq:Xdefinition}-\eqref{eq:Udefinition}, conservative approximations of the chance constraints \eqref{eq:state_and_input_const} are given by the inequality constraints
\begin{subequations}
\begin{align}
    \Pr\left(\alpha_{x,i_x}^\top x_k \leq \beta_{x,i_x} \right) \geq 1 - p_{x,i_x},~~i_x = 1,\ldots,N_s, \label{eq:boole_linear_inequality_state_const}\\
	\Pr\left(\alpha_{u,i_u}^\top u_k \leq \beta_{u,i_u} \right) \geq 1 - p_{u,i_u},~~i_u = 1,\ldots,N_c, \label{eq:boole_linear_inequality_control_const}
\end{align}
\end{subequations}
where $p_{x,i_x}, p_{u,i_u} \ge 0$ are such that
\begin{align}
	\sum_{i_x=1}^{N_s} p_{x,i_x} \leq p_x, \quad
	\sum_{i_u=1}^{N_c} p_{u,i_u} \leq p_u,     
\end{align}
for all $k=0, 1, \dots, N-1$.
We introduce the following lemma to convert \eqref{eq:boole_linear_inequality_state_const} and \eqref{eq:boole_linear_inequality_control_const} to deterministic inequalities. 
\begin{lemma}[\hspace{-0.5pt}\cite{farina2013probabilistic}] \label{lem:cantelli}
    The state chance constraints \eqref{eq:boole_linear_inequality_state_const} are ensured by the tightened deterministic constraints given in terms of the mean and covariance as 
    \begin{align}\label{eq:cantelli_state_const}
        \alpha_{x,i_x}^\top\bar{x}_k + \sqrt{\alpha_{x,i_x}^\top \Sigma_{x_k} \alpha_{x,i_x}}\sqrt{\frac{1 - p_{x, i_x}}{p_{x, i_x}}} - \beta_{x,i_x} \leq 0,
    \end{align}
    where $i_x = 1, \dots, N_s,~~k = 0, \dots, N-1$.
    Likewise, the input constraints \eqref{eq:boole_linear_inequality_control_const} are ensured by the tightened deterministic constraints given in terms of the mean control and control covariance as
    \begin{align}\label{eq:cantelli_input_const}
        \alpha_{u,i_u}^\top\bar{u}_k + \sqrt{\alpha_{u,i_u}^\top \Sigma_{u_k} \alpha_{u,i_u}}\sqrt{\frac{1 - p_{u, i_u}}{p_{u, i_u}}} - \beta_{u,i_u} \leq 0,
    \end{align}
    where $i_u = 1, \dots, N_c,~~k = 0, \dots, N-1$.
    Furthermore, the satisfaction of the original chance constraints \eqref{eq:state_and_input_const} is ensured by the satisfaction of \eqref{eq:cantelli_state_const} and \eqref{eq:cantelli_input_const}.
\end{lemma}

\begin{proof}
    Cantelli's inequality~\cite{marshall1960multivariate} states, for a scalar random variable $y$ with mean $\bar{y}$ and covariance $\Sigma_y$, the inequality,
    \begin{align}
        \Pr(y \geq \bar{y} + \gamma) \leq \frac{\Sigma_y}{\Sigma_y + \gamma^2},
    \end{align}
    holds for every $\gamma \geq 0$, where $\gamma \in \Re$.
    Consider that $\Pr\left(\alpha_{x,i_x}^\top x_k \geq \beta_{x,i_x} \right) \leq p_{x,i_x}$ is equivalent to \eqref{eq:boole_linear_inequality_state_const}, for all $k = 0, \dots, N-1$ and $i_x = 1, \dots, N_s$. 
    Now let $\gamma_{x, i_x, k} = \beta_{x, i_x} - \alpha_{x, i_x}^\top \bar{x}_k$, where $\beta_{x, i_x} \geq 0$ and $\alpha_{x, i_x}^\top \bar{x}_k \leq \beta_{x, i_x}$. Then, $\Pr\left(\alpha_{x,i_x}^\top x_k \geq \alpha_{x, i_x}^\top \bar{x}_k + \gamma_{x, i_x, k} \right) \leq \frac{\alpha_{x,i_x}^\top \Sigma_{x_k} \alpha_{x,i_x}}{\alpha_{x,i_x}^\top \Sigma_{x_k} \alpha_{x,i_x} + \gamma_{x, i_x, k}^2}$. By taking the square root and substituting for $\gamma_{x, i_x, k}$, we arrive at \eqref{eq:cantelli_state_const}. 
    A similar argument may be used for the input constraint to convert \eqref{eq:boole_linear_inequality_control_const} to \eqref{eq:cantelli_input_const}.
    Then, per Boole's inequality~\cite{prekopa1988boole} and the definitions of $\mathcal{X}$ and $\mathcal{U}$ in \eqref{eq:Xdefinition} and \eqref{eq:Udefinition}, respectively, the satisfaction of \eqref{eq:cantelli_state_const} and \eqref{eq:cantelli_input_const} imply satisfaction of \eqref{eq:state_and_input_const}. Specifically, the following inequalities hold
    \begin{align}
        &\Pr(x_k \notin \mathcal{X}) \leq \sum_{i_x=1}^{N_s} \Pr\left(\alpha_{x,i_x}^\top x_k \geq \beta_{x,i_x} \right) \nonumber\\
        &\leq \sum_{i_x = 1}^{N_s} \frac{\alpha_{x,i_x}^\top \Sigma_{x_k} \alpha_{x,i_x}}{\alpha_{x,i_x}^\top \Sigma_{x_k} \alpha_{x,i_x} + \gamma_{x, i_x, k}^2} \leq \sum_{i_x=1}^{N_s} p_{x, i_x} \leq p_x, \\
        &\Pr(u_k \notin \mathcal{U}) \leq \sum_{i_u=1}^{N_c} \Pr\left(\alpha_{u,i_u}^\top u_k \geq \beta_{u,i_u} \right) \nonumber\\
        &\leq \sum_{i_u = 1}^{N_c} \frac{\alpha_{u,i_u}^\top \Sigma_{u_k} \alpha_{u,i_u}}{\alpha_{u,i_u}^\top \Sigma_{u_k} \alpha_{u,i_u} + \gamma_{u, i_u, k}^2} \leq \sum_{i_u=1}^{N_c} p_{u, i_u} \leq p_u
    \end{align}
    for all $k = 0, \dots, N-1$.
\end{proof}

Thus, a deterministic version of Problem~\eqref{prob:cs_prob} is given by
\begin{subequations} \label{prob:determ_cs}
\begin{align}
    \min_{\mathbf{c}, \mathbf{L}} \quad & \bar{J}(\mu_I, \Sigma_I; \mathbf{c}, \mathbf{L}) = \sum_{k=0}^{N-1} \bar{x}_k^\top Q_k \bar{x}_k + \trace(Q_k \Sigma_{x_k}) \nonumber\\
    &\qquad\qquad\qquad + \bar{u}_k^\top R_k \bar{u}_k + \trace(R_k \Sigma_{u_k}) \label{eq:det_cost}\\
    \text{subject to}& \nonumber\\
    \bar{x}_0 &= \mu_I, \label{eq:det_init_mean_const}\\
    \Sigma_{x_0} &= \Sigma_I, \label{eq:det_init_cov_const}\\
    \bar{x}_{k+1} &= \bar{A} \bar{x}_k + \bar{B} \bar{u}_k + \bar{d}, \label{eq:det_mean_dyn_const}\\
    &\hspace{-40pt}\Sigma_{x_{k+1}} = \bar{A} \Sigma_{x_k} \bar{A}^\top + \bar{A} \Sigma_{x_k u_k} \bar{B}^\top + \bar{B} \Sigma_{x_k u_k}^\top \bar{A}^\top + \bar{B} \Sigma_{u_k} \bar{B}^\top \nonumber\\
    &\hspace{-50pt}+ \sum_{j=1}^{m} (\tilde{A}_j \Sigma_{x_k} \tilde{A}_j^\top
    + \tilde{A}_j \Sigma_{x_k u_k} \tilde{B}_j^\top 
    + \tilde{B}_{j} \Sigma_{x_k u_k}^\top \tilde{A}_j^\top
    + \tilde{B}_{j} \Sigma_{u_k} \tilde{B}_j^\top) \nonumber\\
    &\hspace{-40pt}+ \sum_{j=1}^{m} (\tilde{A}_j \bar{x}_k + \tilde{B}_j \bar{u}_k + \tilde{d}_j)(\tilde{A}_j \bar{x}_k + \tilde{B}_j \bar{u}_k + \tilde{d}_j)^\top, \label{eq:det_cov_dyn_const}\\
    \bar{u}_k &= c_k, \\
    \Sigma_{u_k} &= L_k \Sigma_{x_k} L_k^\top, \label{eq:det_input_cov_dyn_const}\\
    \Sigma_{x_k u_k} &= \Sigma_{x_k} L_k^\top, \label{eq:det_state_input_cov_dyn_const}\\
    &\hspace{-12mm}\alpha_{x,i_x}^\top\bar{x}_k + \sqrt{\alpha_{x,i_x}^\top \Sigma_{x_k} \alpha_{x,i_x}}\sqrt{\frac{1 - p_{x, i_x}}{p_{x, i_x}}} - \beta_{x,i_x} \leq 0, \label{eq:det_state_chance_const}\\
    &\hspace{-13mm}\alpha_{u,i_u}^\top\bar{u}_k + \sqrt{\alpha_{u,i_u}^\top \Sigma_{u_k} \alpha_{u,i_u}}\sqrt{\frac{1 - p_{u, i_u}}{p_{u, i_u}}} - \beta_{u,i_u} \leq 0, \label{eq:det_input_chance_const}\\
    \bar{x}_N &= \mu_F, \label{eq:det_terminal_mean_const}\\
    \Sigma_{x_N} &= \Sigma_F \label{eq:det_terminal_cov_const},
\end{align}
\end{subequations}
for $i_x = 1, \dots, N_s$, $i_u = 1, \dots, N_c$, and $k = 0, \dots, N-1$, where $\mathbf{c} = \{c_0, \dots, c_{N-1}\}$ and $\mathbf{L} = \{L_0, \dots, L_{N-1} \}$. We denote the optimal solution of Problem \eqref{prob:determ_cs} as $(\mathbf{c}^\ast, \mathbf{L}^\ast)$ which generates the stochastic control sequence given by 
\begin{align}
    u^\ast_k = L^\ast_k \tilde{x}^\ast_k + c^\ast_k,
\end{align}
where, from \eqref{eq:x_tilde_dynamics} and \eqref{eq:cs_dynamics}, follows that
    $\tilde{x}^\ast_{k+1} = (\bar{A} + \bar{B} L_k^\ast) \tilde{x}^\ast_k 
    + \sum_{j=1}^{m} \tilde{A}_{j} q_{j, k} x^\ast_k 
    + \sum_{j=1}^{m} \tilde{B}_{j} q_{j, k} u^\ast_k 
    + \sum_{j=1}^{m} \tilde{d}_{j} q_{j, k}$, 
    $x^\ast_{k+1} = (\bar{A} + \sum_{j=1}^{m} \tilde{A}_{j} q_{j, k}) x^\ast_k 
    + (\bar{B} + \sum_{j=1}^{m} \tilde{B}_{j} q_{j, k}) u^\ast_k 
    + \bar{d} + \sum_{j=1}^{m} \tilde{d}_{j} q_{j, k}$, 
    $\bar{x}^\ast_{k+1} = \bar{A} \bar{x}_k^\ast + \bar{B} \bar{u}_k^\ast + \bar{d}$,
and $\bar{u}_k^\ast = c_k^\ast$ for all $k = 0, \dots, N-1$, and where $\tilde{x}^\ast_0 = x_0 - \mu_I$, $x_0^\ast = x_0$, and $\bar{x}_0^\ast = \mu_I$.

\begin{theorem}\label{thm:deterministic_original_feasibility}
    Let $\mathbf{c}^\# = \{c_0^\#, \dots, c_{N-1}^\#\}$ and $\mathbf{L}^\# = \{L_0^\#, \dots, L_{N-1}^\# \}$ be a feasible solution of Problem~\eqref{prob:determ_cs}.
    Problem \eqref{prob:determ_cs} is a conservative approximation of Problem \eqref{prob:cs_prob}, such that, any feasible solution of Problem~\eqref{prob:determ_cs} will generate a control sequence given by $\mathbf{u} = \{u_k^\#\}_{k=0}^{N-1} = \{L_k^\# \tilde{x}_k^\# + c_k^\# \}_{k=0}^{N-1}$, which is a feasible solution of Problem~\eqref{prob:cs_prob}. 
    Moreover, by letting $\{u_k^\ast\}_{k=0}^{N-1} = \{L_k^\ast \tilde{x}_k^\ast + c_k^\ast \}_{k=0}^{N-1}$, the optimal solution of Problem \eqref{prob:determ_cs}, in particular, is always a feasible solution of Problem \eqref{prob:cs_prob} satisfying the constraints \eqref{eq:cs_init_mean}-\eqref{eq:cs_terminal_cov_const}.
\end{theorem}

\begin{proof}
    The mean and covariance propagation equations \eqref{eq:det_mean_dyn_const} and \eqref{eq:det_cov_dyn_const} are exact representations for the first two moments of system \eqref{eq:cs_dynamics}, given the control policy \eqref{eq:state_feedback}. Therefore, constraints \eqref{eq:det_init_mean_const}, \eqref{eq:det_init_cov_const}, \eqref{eq:det_terminal_mean_const}, and \eqref{eq:det_terminal_cov_const} are exact representations of \eqref{eq:cs_init_mean}, \eqref{eq:cs_init_cov}, \eqref{eq:cs_terminal_mean_const}, and \eqref{eq:cs_terminal_cov_const}, respectively. For the chance constraints, by Lemma \ref{lem:cantelli}, \eqref{eq:det_state_chance_const} and \eqref{eq:det_input_chance_const} are tightened versions of \eqref{eq:cs_state_chance_const} and \eqref{eq:cs_control_chance_const}, respectively, such that the satisfaction of \eqref{eq:det_state_chance_const} and \eqref{eq:det_input_chance_const} implies satisfaction of \eqref{eq:cs_state_chance_const} and \eqref{eq:cs_control_chance_const}, respectively. Thus, since $\{u_k^\#\}_{k=0}^{N-1}$ and $\{u_k^\ast\}_{k=0}^{N-1}$ necessarily satisfy constraints \eqref{eq:det_init_mean_const}-\eqref{eq:det_terminal_cov_const}, they satisfy constraints \eqref{eq:cs_init_mean}-\eqref{eq:cs_terminal_cov_const}.
\end{proof}

\section{Convex Reformulation} \label{sec:convex_reformulation}

Although Problem~\eqref{prob:determ_cs} is deterministic, it still depends on nonconvex constraints.
First, the propagation of the state covariance \eqref{eq:det_cov_dyn_const} is nonconvex due to the multiplicity of $\bar{x}_k$ and $\bar{u}_k$.
Second, the input and state-input covariance constraints \eqref{eq:det_input_chance_const}-\eqref{eq:det_state_input_cov_dyn_const} are nonconvex owing to the multiplicity of $L_k$ and $\Sigma_{x_k}$.
Third, the state and control constraints \eqref{eq:det_state_chance_const}-\eqref{eq:det_input_chance_const} are nonconvex owing to their nonlinear dependence on $\Sigma_{x_k}$ and $\Sigma_{u_k}$.

We show that the first two issues can be overcome by relaxing the covariance propagation to linear matrix inequality (LMI) constraints, which requires us to relax the terminal covariance constraint \eqref{eq:det_terminal_cov_const} to an inequality constraint, such that we merely ensure the terminal covariance satisfies an upper bound.
The third issue is overcome by further tightening the chance constraints such that they can be written as linear inequality constraints. 
\begin{remark}
    Relaxing the terminal covariance constraint to inequality is not an issue in most applications as, in general, the goal is to bound the covariance rather than to drive it to a specific value, as most applications are concerned with designing a controller to reduce the system's uncertainty rather than increase it.
\end{remark}

\subsubsection{Covariance Propagation Relaxation}

The propagation of the state uncertainty is a nonlinear constraint owing to the multiplicities of $\bar{x}_k \bar{x}_k^\top$, $\bar{x}_k \bar{u}_k^\top$, and $\bar{u}_k \bar{u}_k^\top$. However, this problem can be overcome by relaxing the equality constraint to a LMI.
We introduce a new optimization variable $\bar{\Sigma}_{jk}$ such that
\begin{align} \label{eq:change_of_variables_xbar_ubar}
    \bar{\Sigma}_{jk} \succeq (\tilde{A}_j \bar{x}_k + \tilde{B}_j \bar{u}_k + \tilde{d}_j)(\tilde{A}_j \bar{x}_k + \tilde{B}_j \bar{u}_k + \tilde{d}_j)^\top,
\end{align}
which can be written as a positive semidefinite constraint using the Schur complement
\begin{align}
    \begin{bmatrix}
        \bar{\Sigma}_{jk} & \tilde{A}_j \bar{x}_k + \tilde{B}_j \bar{u}_k + \tilde{d}_j \\
        (\tilde{A}_j \bar{x}_k + \tilde{B}_j \bar{u}_k + \tilde{d}_j)^\top & I
    \end{bmatrix} \succeq 0,
\end{align}
for all $j = 1, \dots, m$ and $k = 0, \dots, N-1$.

Similarly, to handle the nonconvexity of constraints \eqref{eq:det_input_cov_dyn_const} and \eqref{eq:det_state_input_cov_dyn_const}, we utilize the following change of variables previously used by \cite{balci2022exact, benedikter2022convex}. 
Let $\bar{\Sigma}_{ux_k}$ be a new optimization variable such that
\begin{align}\label{eq:change_of_variables_cross_covariance}
    \bar{\Sigma}_{ux_k} = L_k \Sigma_{x_k} = \Sigma_{x_k u_k}^\top,
\end{align}
for all $k = 0, \dots, N-1$.
For $\Sigma_{x_k} \succ 0$, the original control policy can then be recovered as
\begin{align}
    L_k &= \bar{\Sigma}_{ux_k} \Sigma_{x_k}^{-1}.
\end{align}
Using \eqref{eq:change_of_variables_cross_covariance}, $\Sigma_{u_k}$ may be written as 
\begin{align} \label{eq:sigma_u_as_sigma_x_inv}
    \Sigma_{u_k} &= \bar{\Sigma}_{ux_k} \Sigma_{x_k}^{-1} \bar{\Sigma}_{ux_k}^\top.
\end{align}
We then introduce a new optimization variable $\bar{\Sigma}_{u_k}$ and relax \eqref{eq:sigma_u_as_sigma_x_inv} to an inequality given by
\begin{align} \label{eq:change_of_variables_control_cov}
    \bar{\Sigma}_{u_k} \succeq \bar{\Sigma}_{ux_k} \Sigma_{x_k}^{-1} \bar{\Sigma}_{ux_k}^\top = \Sigma_{u_k},
\end{align}
which, using the Schur complement, is given by the positive semidefinite constraint
\begin{align}
    \begin{bmatrix}
        \bar{\Sigma}_{u_k} & \bar{\Sigma}_{ux_k} \\
        \bar{\Sigma}_{ux_k}^\top & \Sigma_{x_k}
    \end{bmatrix} \succeq 0,
\end{align}
for all $k = 0, \dots, N-1$.
We may then write the relaxed covariance dynamics as
\begin{align}\label{eq:change_of_variables_state_covariance}
    \bar{\Sigma}_{x_{k+1}} &= \bar{A} \bar{\Sigma}_{x_k} \bar{A}^\top + \bar{A} \bar{\Sigma}_{ux_k}^\top \bar{B}^\top + \bar{B} \Sigma_{ux_k} \bar{A}^\top + \bar{B} \bar{\Sigma}_{u_k} \bar{B}^\top \nonumber\\
        &+ \sum_{j=1}^{m} (\tilde{A}_j \bar{\Sigma}_{x_k} \tilde{A}_j^\top
        + \tilde{A}_j \Sigma_{ux_k}^\top \tilde{B}_j^\top 
        + \tilde{B}_{j} \Sigma_{ux_k} \tilde{A}_j^\top \nonumber\\
        &+ \tilde{B}_{j} \bar{\Sigma}_{u_k} \tilde{B}_j^\top 
        + \bar{\Sigma}_{jk}),
\end{align}
for $k = 0, \dots, N-1$, where $\bar{\Sigma}_{x_0} = \Sigma_{x_0}$.

Relaxing \eqref{eq:det_cov_dyn_const} and \eqref{eq:det_input_cov_dyn_const} with \eqref{eq:change_of_variables_xbar_ubar}, \eqref{eq:change_of_variables_cross_covariance}, \eqref{eq:change_of_variables_control_cov}, and \eqref{eq:change_of_variables_state_covariance}
provides a bound on the state and input covariances. 
\begin{lemma}\label{lem:cov_dynamics_relaxation}
    If the inequalities 
    \begin{subequations}
    \begin{align}
        \alpha_{x,i_x}^\top\bar{x}_k &+ \sqrt{\alpha_{x,i_x}^\top \bar{\Sigma}_{x_k} \alpha_{x,i_x}}\sqrt{\frac{1 - p_{x, i_x}}{p_{x, i_x}}} - \beta_{x,i_x} \leq 0, \label{eq:new_state_chance_const_upper_bounds} \\
        \alpha_{u,i_u}^\top\bar{u}_k &+ \sqrt{\alpha_{u,i_u}^\top \bar{\Sigma}_{u_k} \alpha_{u,i_u}}\sqrt{\frac{1 - p_{u, i_u}}{p_{u, i_u}}} - \beta_{u,i_u} \leq 0,  \label{eq:new_input_chance_const_upper_bound} \\
        \bar{\Sigma}_{x_N} &\preceq \Sigma_{F}. \label{eq:new_terminal_cov_upper_bound}
    \end{align}
    \end{subequations}
    hold, 
    then the inequalities 
    \begin{subequations}
    \begin{align} 
        \alpha_{x,i_x}^\top\bar{x}_k &+ \sqrt{\alpha_{x,i_x}^\top \Sigma_{x_k} \alpha_{x,i_x}}\sqrt{\frac{1 - p_{x, i_x}}{p_{x, i_x}}} - \beta_{x,i_x} \leq 0, \label{eq:old_state_chance_const_upper_bounds} \\
        \alpha_{u,i_u}^\top\bar{u}_k &+ \sqrt{\alpha_{u,i_u}^\top \Sigma_{u_k} \alpha_{u,i_u}}\sqrt{\frac{1 - p_{u, i_u}}{p_{u, i_u}}} - \beta_{u,i_u} \leq 0, \label{eq:old_input_chance_const_upper_bound} \\
        \Sigma_{x_N} &\preceq \Sigma_{F}. \label{eq:old_terminal_cov_upper_bound}
    \end{align}
    \end{subequations}
    also hold for $i_x = 1, \dots, N_s$, $i_u = 1, \dots, N_c$, and $k = 0, \dots, N-1$.
\end{lemma}
\begin{proof}
    Given \eqref{eq:change_of_variables_xbar_ubar}, \eqref{eq:change_of_variables_cross_covariance}, and \eqref{eq:change_of_variables_control_cov}, it can be seen that the following inequality holds,
    \begin{align}\label{eq:psd_covariance}
        \bar{\Sigma}_{x_{k+1}} &= \bar{A} \bar{\Sigma}_{x_k} \bar{A}^\top + \bar{A} \bar{\Sigma}_{ux_k}^\top \bar{B}^\top + \bar{B} \bar{\Sigma}_{ux_k} \bar{A}^\top + \bar{B} \bar{\Sigma}_{u_k} \bar{B}^\top \nonumber\\
        &+ \sum_{j=1}^{m} (\tilde{A}_j \bar{\Sigma}_{x_k} \tilde{A}_j^\top
        + \tilde{A}_j \bar{\Sigma}_{ux_k}^\top \tilde{B}_j^\top 
        + \tilde{B}_{j} \bar{\Sigma}_{ux_k} \tilde{A}_j^\top \nonumber\\
        &+ \tilde{B}_{j} \bar{\Sigma}_{u_k} \tilde{B}_j^\top 
        + \bar{\Sigma}_{jk}) \nonumber\\
        \succeq \Sigma_{x_{k+1}} &= \bar{A} \Sigma_{x_k} \bar{A}^\top + \bar{A} \Sigma_{x_k u_k} \bar{B}^\top + \bar{B} \Sigma_{x_k u_k}^\top \bar{A}^\top + \bar{B} \Sigma_{u_k} \bar{B}^\top \nonumber\\
        &+ \sum_{j=1}^{m} (\tilde{A}_j \Sigma_{x_k} \tilde{A}_j^\top
        + \tilde{A}_j \Sigma_{x_k u_k} \tilde{B}_j^\top 
        + \tilde{B}_{j} \Sigma_{x_k u_k}^\top \tilde{A}_j^\top  \nonumber\\
        &+ \tilde{B}_{j} \Sigma_{u_k} \tilde{B}_j^\top)
        + \sum_{j=1}^{m} (\tilde{A}_j \bar{x}_k + \tilde{B}_j \bar{u}_k + \tilde{d}_j)(\tilde{A}_j \bar{x}_k  \nonumber\\
        &+ \tilde{B}_j \bar{u}_k + \tilde{d}_j)^\top,
    \end{align}
    for $k = 0, \dots, N-1$.
    Therefore, given \eqref{eq:new_state_chance_const_upper_bounds}, \eqref{eq:new_input_chance_const_upper_bound}, and \eqref{eq:change_of_variables_control_cov}, and since $\Sigma_{x_k} \succ 0$, $\Sigma_{u_k} \succeq 0$, and $\sqrt{\frac{1 - p_{x, i_x}}{p_{x, i_x}}} \geq 0$, the following inequalities hold,
    \begin{subequations}
    \begin{align}
        &\alpha_{x,i_x}^\top\bar{x}_k + \sqrt{\alpha_{x,i_x}^\top \Sigma_{x_k} \alpha_{x,i_x}}\sqrt{\frac{1 - p_{x, i_x}}{p_{x, i_x}}} - \beta_{x,i_x} \nonumber\\
        &\quad \leq \alpha_{x,i_x}^\top\bar{x}_k + \sqrt{\alpha_{x,i_x}^\top \bar{\Sigma}_{x_k} \alpha_{x,i_x}}\sqrt{\frac{1 - p_{x, i_x}}{p_{x, i_x}}} - \beta_{x,i_x} \leq 0, \\
        &\alpha_{u,i_u}^\top\bar{u}_k + \sqrt{\alpha_{u,i_u}^\top \Sigma_{u_k} \alpha_{u,i_u}}\sqrt{\frac{1 - p_{u, i_u}}{p_{u, i_u}}} - \beta_{u,i_u} \nonumber\\
        &\quad \leq  \alpha_{u,i_u}^\top\bar{u}_k + \sqrt{\alpha_{u,i_u}^\top \bar{\Sigma}_{u_k} \alpha_{u,i_u}}\sqrt{\frac{1 - p_{u, i_u}}{p_{u, i_u}}} - \beta_{u,i_u} \leq 0,
    \end{align}
    \end{subequations}
    for $i_x = 1, \dots, N_s$, $ i_u = 1, \dots, N_c$, and $k = 0, \dots, N-1$.
    Additionally, if $\Sigma_F \succeq \bar{\Sigma}_{x_N}$, then the inequality $\Sigma_F \succeq \bar{\Sigma}_{x_N} \succeq \Sigma_{x_N}$ holds.
\end{proof}

\subsubsection{Chance Constraint Tightening}

Finally, the state and input chance constraints \eqref{eq:det_state_chance_const}-\eqref{eq:det_input_chance_const} include a nonlinear dependence on the state and input covariance, $\Sigma_{x_k}$ and $\Sigma_{u_k}$, respectively. To resolve this, we use a linearization procedure similar to \cite{farina2013probabilistic}. An upper bound on the square-root of the state covariance can be derived using a tangent line approximation evaluated at $\lambda_{x, i_x, k}$ given by
\begin{align}\label{eq:sqrt_state_cov_lin_approx}
    \sqrt{\alpha_{x,i_x}^\top \Sigma_{x_k} \alpha_{x,i_x}} &\leq \frac{1}{2 \sqrt{\lambda_{x, i_x, k}}} (\alpha_{x,i_x}^\top \Sigma_{x_k} \alpha_{x,i_x} - \lambda_{x, i_x, k}) \nonumber\\
    &+ \sqrt{\lambda_{x, i_x, k}},
\end{align}
for $i_x = 1, \dots, N_s$ and $k = 0, \dots, N-1$.
A similar bound is given for the input covariance as 
\begin{align}\label{eq:sqrt_input_cov_lin_approx}
    \sqrt{\alpha_{u,i_u}^\top \Sigma_{u_k} \alpha_{u,i_u}} &\leq \frac{1}{2 \sqrt{\lambda_{u, i_u, k}}} (\alpha_{u,i_u}^\top \Sigma_{u_k} \alpha_{u,i_u} - \lambda_{u, i_u, k}) \nonumber\\
    &+ \sqrt{\lambda_{u, i_u, k}},
\end{align}
for $i_u = 1, \dots, N_c$ and $k=0, \dots, N-1$.
This leads to a conservative tightening of the chance constraints; therefore, the original constraints are still guaranteed to be satisfied by a solution of the convex reformulation at the expense of reducing the size of the feasible solution set, as stated in the following lemma.
\begin{lemma}\label{lem:chance_const_linearization}
    The satisfaction of the inequalities given by
    \begin{subequations}
    \begin{align}
        &\alpha_{x,i_x}^\top\bar{x}_k + (\frac{\sqrt{\lambda_{x, i_x, k}}}{2} + \frac{1}{2\sqrt{\lambda_{x, i_x, k}}} \alpha_{x,i_x}^\top \bar{\Sigma}_{x_k} \alpha_{x,i_x})\sqrt{\frac{1 - p_{x, i_x}}{p_{x, i_x}}} \nonumber\\
        &- \beta_{x,i_x} \leq 0, \label{eq:linearized_state_chance}\\
        &\alpha_{u,i_u}^\top\bar{u}_k + (\frac{\sqrt{\lambda_{u, i_u, k}}}{2} + \frac{1}{2\sqrt{\lambda_{u, i_u, k}}} \alpha_{u,i_u}^\top \bar{\Sigma}_{u_k} \alpha_{u,i_u}) \sqrt{\frac{1 - p_{u, i_u}}{p_{u, i_u}}} \nonumber\\
        &- \beta_{u,i_u} \leq 0, \label{eq:linearized_input_chance}
    \end{align} 
    \end{subequations}
    for $i_x = 1, \dots, N_s$, $i_u = 1, \dots, N_c$, and $k = 0, \dots, N-1$, imply satisfaction of \eqref{eq:det_state_chance_const} and \eqref{eq:det_input_chance_const}.   
\end{lemma}
\begin{proof}
    Considering the concavity of $\sqrt{\alpha_{x,i_x}^\top \bar{\Sigma}_{x_k} \alpha_{x,i_x}}$ and $\sqrt{\alpha_{u,i_u}^\top \bar{\Sigma}_{u_k} \alpha_{u,i_u}}$ and that $\bar{\Sigma}_{x_k} \succ 0$ and $\bar{\Sigma}_{u_k} \succeq 0$, it follows that \eqref{eq:sqrt_state_cov_lin_approx} and \eqref{eq:sqrt_input_cov_lin_approx} hold. Moreover, since $\sqrt{\frac{1 - p_{x, i_x}}{p_{x, i_x}}} \geq 0$, if \eqref{eq:linearized_state_chance} and \eqref{eq:linearized_input_chance} hold, then the inequalities
    \begin{subequations}
    \begin{align}
        &\alpha_{x,i_x}^\top\bar{x}_k - \beta_{x,i_x}
        + \sqrt{\alpha_{x,i_x}^\top \bar{\Sigma}_{x_k} \alpha_{x,i_x}}\sqrt{\frac{1 - p_{x, i_x}}{p_{x, i_x}}} \nonumber\\
        &\leq (\frac{\sqrt{\lambda_{x, i_x, k}}}{2} + \frac{1}{2\sqrt{\lambda_{x, i_x, k}}} \alpha_{x,i_x}^\top \bar{\Sigma}_{x_k} \alpha_{x,i_x})\sqrt{\frac{1 - p_{x, i_x}}{p_{x, i_x}}} \nonumber\\
        &+ \alpha_{x,i_x}^\top\bar{x}_k - \beta_{x,i_x} \leq 0, \label{eq:original_and_linearized_state_chance}\\
        &\alpha_{u,i_u}^\top\bar{u}_k + \sqrt{\alpha_{u,i_u}^\top \bar{\Sigma}_{u_k} \alpha_{u,i_u}}\sqrt{\frac{1 - p_{u, i_u}}{p_{u, i_u}}} - \beta_{u,i_u} \nonumber\\
        &\leq (\frac{\sqrt{\lambda_{u, i_u, k}}}{2} + \frac{1}{2\sqrt{\lambda_{u, i_u, k}}} \alpha_{u,i_u}^\top \bar{\Sigma}_{u_k} \alpha_{u,i_u}) \sqrt{\frac{1 - p_{u, i_u}}{p_{u, i_u}}} \nonumber\\
        &+ \alpha_{u,i_u}^\top\bar{u}_k - \beta_{u,i_u} \leq 0 \label{eq:original_and_linearized_input_chance}
    \end{align}
    \end{subequations}
    hold for all $i_x = 1, \dots, N_s$, $i_u = 1, \dots, N_c$, and $k = 0, \dots, N-1$, which by Lemma~\ref{lem:cov_dynamics_relaxation} imply satisfaction of \eqref{eq:det_state_chance_const} and \eqref{eq:det_input_chance_const}.
\end{proof}

\begin{remark}
    The approximations given by \eqref{eq:sqrt_state_cov_lin_approx} and \eqref{eq:sqrt_input_cov_lin_approx} are exact (i.e., satisfied as an equality) when $\lambda_{x, i_x, k} = \alpha_{x,i_x}^\top \Sigma_{x_k} \alpha_{x,i_x}$ or $\lambda_{u, i_u, k} = \alpha_{u,i_u}^\top \Sigma_{u_k} \alpha_{u,i_u}$, respectively. Therefore, if estimates of $\Sigma_{x_k}$ and $\Sigma_{u_k}$ are available, they should be used for the selection of the linearization points (e.g., by interpolating between the given initial and final covariances or solving the problem iteratively until convergence, where the values from the prior solution are used as the linearization points).
\end{remark}

\subsubsection{Convex Covariance Steering Problem}

Using the above relaxations, a convex formulation of Problem~\eqref{prob:determ_cs} is given by
\begin{subequations} \label{prob:convex_cs}
\begin{align}
    \min_{\mathbf{c}, \mathbf{\bar{\Sigma}_{ux}}, \mathbf{\bar{\Sigma}_u}, \mathbf{\bar{\Sigma}}} & \bar{x}_k^\top Q_k \bar{x}_k + \trace(Q_k \bar{\Sigma}_{x_k}) + \bar{u}_k^\top R_k \bar{u}_k + \trace(R_k \bar{\Sigma}_{u_k}) \\
    \text{subject to}& \nonumber\\
    \bar{x}_0 &= \mu_I, \label{eq:convex_init_mean_const}\\
    \bar{\Sigma}_{x_0} &= \Sigma_I, \label{eq:convex_init_cov_const}\\
    \bar{x}_{k+1} &= \bar{A} \bar{x}_k + \bar{B} \bar{u}_k + \bar{d} \label{eq:convex_mean_dynamics_const}\\
    \bar{\Sigma}_{x_{k+1}} &= \bar{A} \bar{\Sigma}_{x_k} \bar{A}^\top + \bar{A} \bar{\Sigma}_{ux_k}^\top \bar{B}^\top + \bar{B} \bar{\Sigma}_{ux_k} \bar{A}^\top \nonumber\\
    &+ \bar{B} \bar{\Sigma}_{u_k} \bar{B}^\top
    + \sum_{j=1}^{m} (\tilde{A}_j \bar{\Sigma}_{x_k} \tilde{A}_j^\top
    + \tilde{A}_j \bar{\Sigma}_{ux_k}^\top \tilde{B}_j^\top \nonumber\\
    &+ \tilde{B}_{j} \bar{\Sigma}_{ux_k} \tilde{A}_j^\top
    + \tilde{B}_{j} \bar{\Sigma}_{u_k} \tilde{B}_j^\top
    + \bar{\Sigma}_{jk}), \label{eq:convex_cov_dynamics_const}\\
    \bar{u}_k &= c_k \label{eq:convex_ubar_const}\\
    \bar{\Sigma}_{u_k} &\succeq \bar{\Sigma}_{ux_k} \bar{\Sigma}_{x_k}^{-1} \bar{\Sigma}_{ux_k}^\top \label{eq:convex_sigma_u_const}\\
    &\hspace{-35pt} \bar{\Sigma}_{jk} \succeq (\tilde{A}_j \bar{x}_k + \tilde{B}_j \bar{u}_k + \tilde{d}_j)(\tilde{A}_j \bar{x}_k + \tilde{B}_j \bar{u}_k + \tilde{d}_j)^\top \label{eq:convex_sigma_bar_const}\\
    &\hspace{-15pt}(\frac{\sqrt{\lambda_{x, i_x, k}}}{2} + \frac{1}{2\sqrt{\lambda_{x, i_x, k}}} \alpha_{x,i_x}^\top \bar{\Sigma}_{x_k} \alpha_{x,i_x})\sqrt{\frac{1 - p_{x, i_x}}{p_{x, i_x}}} \nonumber\\
    &+ \alpha_{x,i_x}^\top\bar{x}_k - \beta_{x,i_x} \leq 0, \label{eq:convex_state_chance_const}\\
    &\hspace{-15pt}(\frac{\sqrt{\lambda_{u, i_u, k}}}{2} + \frac{1}{2\sqrt{\lambda_{u, i_u, k}}} \alpha_{u,i_u}^\top \bar{\Sigma}_{u_k} \alpha_{u,i_u}) \sqrt{\frac{1 - p_{u, i_u}}{p_{u, i_u}}} \nonumber\\
    &+ \alpha_{u,i_u}^\top\bar{u}_k - \beta_{u,i_u} \leq 0, \label{eq:convex_input_chance_const}\\
    \bar{x}_N &= \mu_F, \label{eq:convex_terminal_mean_const}\\
    \bar{\Sigma}_{x_N} &\preceq \Sigma_F, \label{eq:convex_terminal_cov_const}
\end{align}
\end{subequations}
for $i_x = 1, \dots, N_s$, $i_u = 1, \dots, N_u$, and $k=0, \dots, N-1$, where $\mathbf{\bar{\Sigma}_{ux}} = \{ \bar{\Sigma}_{ux_0}, \dots, \bar{\Sigma}_{ux_{N-1}} \}$, $\mathbf{\bar{\Sigma}_u} = \{ \bar{\Sigma}_{u_0}, \dots, \bar{\Sigma}_{u_{N-1}} \}$, and $\mathbf{\bar{\Sigma}} = \{\bar{\Sigma}_{j, k} \}_{j=1, k=0}^{m, N-1}$.

Let us denote the optimal solution of Problem \eqref{prob:convex_cs} as $(\mathbf{c}^\star, \mathbf{\bar{\Sigma}_{ux}}^\star, \mathbf{\bar{\Sigma}_u}^\star, \mathbf{\bar{\Sigma}}^\star)$, which gives the optimal control policy $(\mathbf{c}^\star, \mathbf{L}^\star)$, where 
    $L_k^\star = \bar{\Sigma}_{ux_k}^\star \bar{\Sigma}_{x_k}^{\star-1}$, 
    $\bar{\Sigma}_{x_{k+1}}^\star = \bar{A} \bar{\Sigma}_{x_k}^\star \bar{A}^\top + \bar{A} \bar{\Sigma}_{ux_k}^{\star\top} \bar{B}^\top + \bar{B} \bar{\Sigma}_{ux_k}^\star \bar{A}^\top + \bar{B} \bar{\Sigma}_{u_k}^\star \bar{B}^\top 
    + \sum_{j=1}^{m} (\tilde{A}_j \bar{\Sigma}_{x_k}^\star \tilde{A}_j^\top
    + \tilde{A}_j \bar{\Sigma}_{ux_k}^{\star\top} \tilde{B}_j^\top 
    + \tilde{B}_{j} \bar{\Sigma}_{ux_k}^\star \tilde{A}_j^\top 
    + \tilde{B}_{j} \bar{\Sigma}_{u_k}^\star \tilde{B}_j^\top
    + \bar{\Sigma}_{jk}^\star)$, 
and $\bar{\Sigma}_{x_0}^\star = \Sigma_I$.

\begin{theorem}
    The optimal solution of Problem~\eqref{prob:convex_cs}, $(\mathbf{c}^\star, \mathbf{\bar{\Sigma}_{ux}}^\star, \mathbf{\bar{\Sigma}_u}^\star, \mathbf{\bar{\Sigma}}^\star)$, yields the optimal control policy $\{c_k^\star, L_k^\star\}_{k=0}^{N-1}$ which is a feasible solution of Problem \eqref{prob:determ_cs} when constraint \eqref{eq:det_terminal_cov_const} is relaxed to \eqref{eq:old_terminal_cov_upper_bound}. Furthermore, let $\mathbf{u}^\star$ be the control sequence given by $\mathbf{u}^\star = \{u_k^\star = L_k^\star \tilde{x}_k^\star + c_k^\star \}_{k=0}^{N-1}$,
    where 
        $\tilde{x}^\star_{k+1} = (\bar{A} + \bar{B} L_k^\star) \tilde{x}^\star_k 
        + \sum_{j=1}^{m} \tilde{A}_{j} q_{j, k} x^\star_k 
        + \sum_{j=1}^{m} \tilde{B}_{j} q_{j, k} u^\star_k
        + \sum_{j=1}^{m} \tilde{d}_{j} q_{j, k}$, 
        $x^\star_{k+1} = (\bar{A} + \sum_{j=1}^{m} \tilde{A}_{j} q_{j, k}) x^\star_k + (\bar{B} + \sum_{j=1}^{m} \tilde{B}_{j} q_{j, k}) u^\star_k + \bar{d} + \sum_{j=1}^{m} \tilde{d}_{j} q_{j, k}$, 
        $\bar{x}^\star_{k+1} = \bar{A} \bar{x}_k^\star + \bar{B} \bar{u}_k^\star + \bar{d}$,
    and $\bar{u}_k^\star = c_k^\star$ for all $k = 0, \dots, N-1$, and where $\tilde{x}^\star_0 = x_0 - \mu_I$, $x_0^\star = x_0$, and $\bar{x}_0^\star = \mu_I$.
    The control sequence $\mathbf{u}^\star$ is a feasible solution of Problem \eqref{prob:cs_prob} when \eqref{eq:cs_terminal_cov_const} is relaxed to \eqref{eq:old_terminal_cov_upper_bound}.
\end{theorem}
\begin{proof}
    Constraints \eqref{eq:convex_init_mean_const}, \eqref{eq:convex_init_cov_const}, \eqref{eq:convex_mean_dynamics_const}, \eqref{eq:convex_ubar_const}, and \eqref{eq:convex_terminal_mean_const} are identical to the corresponding constraints in Problem \eqref{prob:determ_cs}. 
    %
    Per Lemma \ref{lem:cov_dynamics_relaxation}, satisfaction of \eqref{eq:convex_terminal_cov_const} implies \eqref{eq:old_terminal_cov_upper_bound} will be satisfied.
    By Lemmas \ref{lem:cov_dynamics_relaxation} and \ref{lem:chance_const_linearization} satisfaction of \eqref{eq:convex_state_chance_const} and \eqref{eq:convex_input_chance_const} implies satisfaction of \eqref{eq:det_state_chance_const} and \eqref{eq:det_input_chance_const}.
    Thus, $\{c_k^\star, L_k^\star\}_{k=0}^{N-1}$ is a feasible solution of Problem~\eqref{prob:determ_cs} when constraint \eqref{eq:det_terminal_cov_const} is relaxed to \eqref{eq:old_terminal_cov_upper_bound}.
    Feasibility of Problem~\eqref{prob:cs_prob} with the relaxation of \eqref{eq:cs_terminal_cov_const} to \eqref{eq:old_terminal_cov_upper_bound} is then given by Theorem~\ref{thm:deterministic_original_feasibility}.
\end{proof}

\section{Numerical Results} \label{sec:numerical_results}

\begin{figure}[h]
    \centering
    \includegraphics[trim= 40 80 40 120,clip,width=\columnwidth]{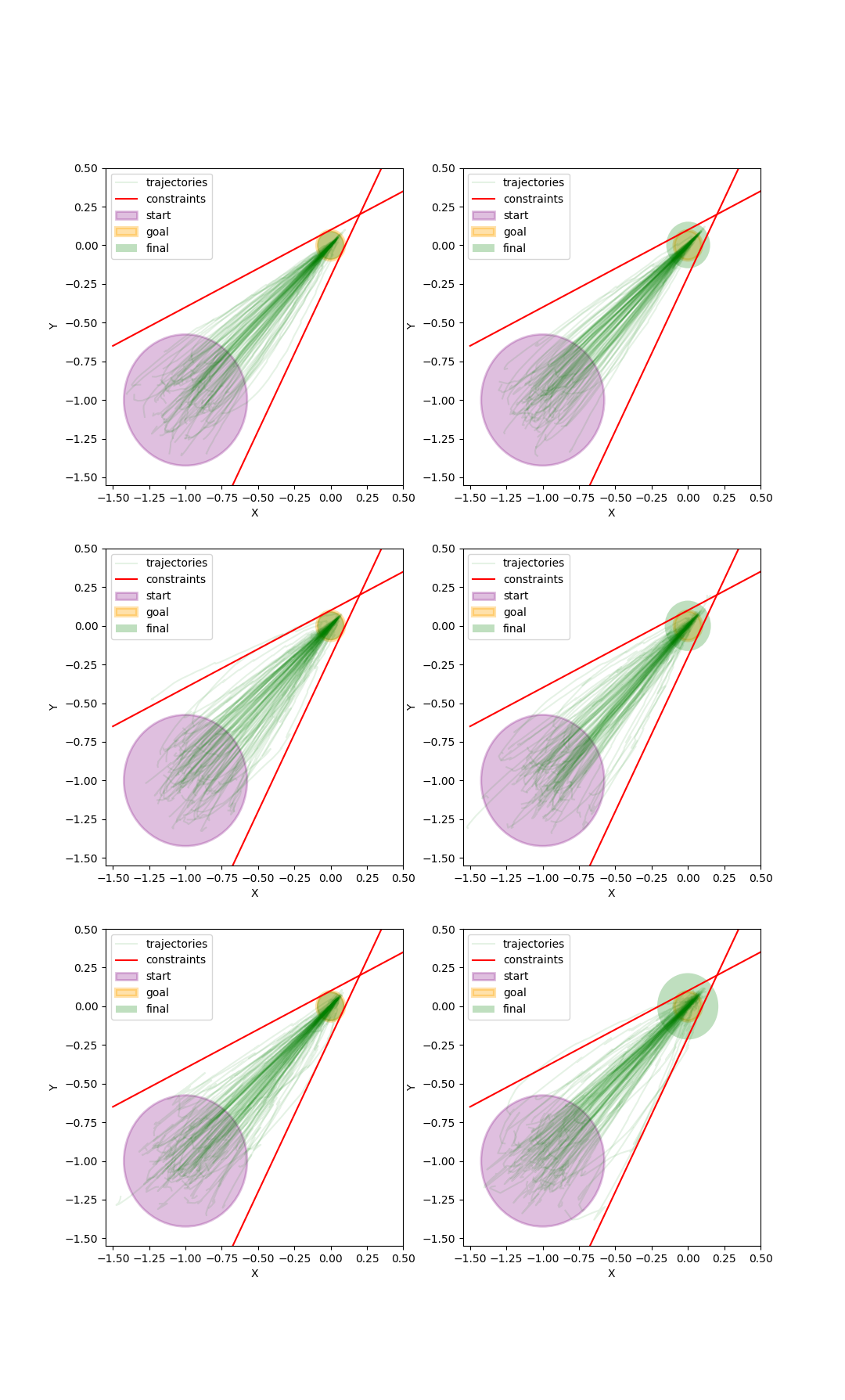}
    \caption{Monte Carlo trajectory simulations for different noise intensities (rows) and problem formulations (columns). The left column shows the proposed approach defined in Problem \eqref{prob:convex_cs}, whereas the right column shows a naive approach that does not consider the multiplicative uncertainties arising from $\tilde{A}_j$ and $\tilde{B}_j$. The rows show the results for noise intensities $\theta = 0.1, 0.5, 1.0$, respectively. It may be seen that the proposed approach satisfies the terminal covariance constraint whereas the naive approach does not.}
    \label{fig:monte_carlo_simulations}
\end{figure}

The convex problem formulation \eqref{prob:convex_cs} is verified on a 2D double integrator system, and the results are compared with a naive problem formulation which neglects to consider the multiplicative uncertainties. The double integrator system is given by 
\begin{subequations}
\begin{align*}
    \bar{A} = \begin{bmatrix}
        1 & 0 & 0 & 0 \\
        0 & 1 & 0 &0 \\
        \Delta t & 0 & 1 & 0 \\
        0 & \Delta t & 0 & 1
    \end{bmatrix},
    & \quad
    \tilde{A}_1 = \begin{bmatrix}
        0 & 0 & 0 & 0 \\
        0 & 0 & 0 & 0 \\
        \theta \Delta t & 0 & 0 & 0 \\
        0 & \theta \Delta t & 0 & 0 \\
    \end{bmatrix},
    \\
    \tilde{A}_2 = \begin{bmatrix}
        \theta \Delta t & 0 & 0 & 0 \\
        0 & \theta \Delta t & 0 & 0 \\
        0 & 0 & 0 & 0 \\
        0 & 0 & 0 & 0 \\
    \end{bmatrix},
    & \quad
    \bar{B} = \begin{bmatrix}
         \Delta t & 0 \\
         0 & \Delta t \\
         0 & 0 \\
         0 & 0
    \end{bmatrix},
    \\
    \tilde{B}_1 = \begin{bmatrix}
         \theta \Delta t & 0 \\
         0 & \theta \Delta t \\
         0 & 0 \\
         0 & 0
    \end{bmatrix},
    & \quad
    \tilde{B}_2 = \begin{bmatrix}
         0 & \theta \Delta t \\
         \theta \Delta t & 0 \\
         0 & 0 \\
         0 & 0
    \end{bmatrix},
    \\
    \bar{d} = 0_4,
    \quad 
    \tilde{d}_1 = \begin{bmatrix}
        0.1 \Delta t \\
        0.1 \Delta t \\
        0 \\
        0
    \end{bmatrix},
    & \quad 
    \tilde{d}_1 = \begin{bmatrix}
        0 \\
        0 \\
        0.1 \Delta t \\
        0.1 \Delta t \\
    \end{bmatrix},
\end{align*}
\end{subequations}
where $\Delta t = 0.1 \text{s}$ is the time-step, and $\theta$ is the noise intensity which is varied between simulations and given in Figure \ref{fig:monte_carlo_simulations}. The initial condition is given by $\mu_I = [0, 0, -1, -1]^\top$, $\Sigma_I = 0.02 I_4$ and the terminal constraints are given as $\mu_F = [0, 0, 0, 0]^\top$, $\Sigma_F = 10 (\tilde{d}_1 \tilde{d}_1^\top + \tilde{d}_2 \tilde{d}_2^\top)$. The state chance constraints are given by $\alpha_{x, 1} = [0, 0, -0.5, 1.0]^\top$, $\beta_{x, 1} = 0.1$, $p_{x, 1} = 0.2$, $\alpha_{x, 2} = [0, 0, 2.0, -1.0]^\top$, $\beta_{x, 2} = 0.2$, $p_{x, 2} = 0.2$ and no control chance constraints are included so that $N_s = 2$ and $N_c = 0$. $\lambda_{x, i_x, k} = \alpha_{x, i_x}^\top \Sigma^{\text{nom}}_{x,k} \alpha_{x, i_x}$ for $i_x = 1, 2$ and $k = 0, \dots, N-1$, and where $\Sigma^{\text{nom}}_{x,k}$ is computed as an $N$-step linear interpolation between $\Sigma_I$ and $\Sigma_F$. The trajectory is planned over 5 seconds, so that $N = 50$.

\section{Conclusion} \label{sec:conclusion}

This work has presented a general problem formulation for stochastic linear systems subjected to both additive and multiplicative noise, while subject to state and control chance constraints as well as terminal constraints on the first and second moments. Although the problem is, in general, stochastic and nonconvex, a tightened, deterministic, convex problem formulation is derived, the optimal solution of which is guaranteed to be a valid (albeit potentially sub-optimal) solution of the original nonconvex problem. Finally, the results are demonstrated using Monte Carlo simulations on a 2D double integrator system. 

\bibliographystyle{IEEEtran}
\bibliography{parametric_uncertainty.bib}

\end{document}